\documentclass[a4paper,12pt]{amsart}

\usepackage{amsmath}
\usepackage{amssymb}
\usepackage{mathrsfs}
\usepackage{ifthen}
\usepackage{graphicx}
\usepackage[T1]{fontenc} 

\def\switchlinenumbers{\@ifstar
    {\let\makeLineNumberOdd\makeLineNumberRight
     \let\makeLineNumberEven\makeLineNumberLeft}%
    {\let\makeLineNumberOdd\makeLineNumberLeft
     \let\makeLineNumberEven\makeLineNumberRight}%
    }

\def\setmakelinenumbers#1{\@ifstar
  {\let\makeLineNumberRunning#1%
   \let\makeLineNumberOdd#1%
   \let\makeLineNumberEven#1}%
  {\ifx\c@linenumber\c@runninglinenumber
      \let\makeLineNumberRunning#1%
   \else
      \let\makeLineNumberOdd#1%
      \let\makeLineNumberEven#1%
   \fi}%
  }

\nonstopmode \numberwithin{equation}{section}
\newtheorem*{theorem*}{Theorem}

\newtheorem{thm}{Theorem}[section]
\newtheorem{cor}{Corollary}[section]
\newtheorem{lem}{Lemma}[section]

\theoremstyle{definition}
\newtheorem{defn}{Definition}[section]

\newtheorem{prob}[equation]{Problem}
\newtheorem{rem}{Remark}[section]

\newtheorem{innercustomthm}{Theorem}
\newenvironment{customthm}[1]
  {\renewcommand\theinnercustomthm{#1}\innercustomthm}
  {\endinnercustomthm}

\newcounter{minutes}
\divide\time by 60
\newcounter{hours}
\multiply\time by 60
\addtocounter{minutes}{-\time}

\newcounter {own}
\def\theown {\thesection       .\arabic{own}}

\newenvironment{pf}[1][]{%
 \vskip 3mm
 \noindent
 \ifthenelse{\equal{#1}{}}%
  {{\slshape Proof. }}%
  {{\slshape #1.} }%
 }%
{\qed\bigskip}

\newcounter{alphabet}

\def\be{\begin{equation}}
\def\ee{\end{equation}}

\newcommand{\bee}{\begin{enumerate}}
\newcommand{\eee}{\end{enumerate}}

\newcommand{\blem}{\begin{lem}}
\newcommand{\elem}{\end{lem}}
\newcommand{\bthm}{\begin{thm}}
\newcommand{\ethm}{\end{thm}}
\newcommand{\bcor}{\begin{cor}}
\newcommand{\ecor}{\end{cor}}
\newcommand{\beg}{\begin{examp}}
\newcommand{\eeg}{\end{examp}}
\newcommand{\begs}{\begin{examples}}
\newcommand{\eegs}{\end{examples}}
\newcommand{\bdefe}{\begin{defin}}
\newcommand{\edefe}{\end{defin}}
\newcommand{\bprob}{\begin{prob}}
\newcommand{\eprob}{\end{prob}}
\newcommand{\bei}{\begin{itemize}}
\newcommand{\eei}{\end{itemize}}

\newcommand{\norm}[1]{\left\lVert#1\right\rVert}

\newcommand{\innpdct}[1]{\left\langle#1\right\rangle}

\begin{document}

\title{Bohr and Rogosinski inequalities for operator valued holomorphic functions}

\author{Vasudevarao Allu}
\address{Vasudevarao Allu,
School of Basic Sciences,
Indian Institute of Technology Bhubaneswar,
Bhubaneswar-752050, Odisha, India.}
\email{avrao@iitbbs.ac.in}

\author{Himadri Halder}
\address{Himadri Halder,
School of Basic Sciences,
Indian Institute of Technology Bhubaneswar,
Bhubaneswar-752050, Odisha, India.}
\email{himadrihalder119@gmail.com}

\author{Subhadip Pal}
\address{Subhadip Pal,
	School of Basic Sciences,
	Indian Institute of Technology Bhubaneswar,
	Bhubaneswar-752050, Odisha, India.}
\email{subhadippal33@gmail.com}

\subjclass[{AMS} Subject Classification:]{Primary 46E40, 47A56; Secondary 46B20, 47A63, 30B10}

\def\thefootnote{}
\footnotetext{ {\tiny File:~\jobname.tex,
printed: \number\year-\number\month-\number\day,
          \thehours.\ifnum\theminutes<10{0}\fi\theminutes }
} \makeatletter\def\thefootnote{\@arabic\c@footnote}\makeatother

\begin{abstract}
 For any complex Banach space $X$ and each $p \in [1,\infty)$, we introduce the $p$-Bohr radius of order $N(\in \mathbb{N})$ is $\widetilde{R}_{p,N}(X)$ defined by
$$
\widetilde{R}_{p,N}(X)=\sup \left\{r\geq 0: \sum_{k=0}^{N}\norm{x_k}^p r^{pk} \leq \norm{f}^p_{H^{\infty}(\mathbb{D}, X)}\right\},
$$
where $f(z)=\sum_{k=0}^{\infty} x_{k}z^k \in H^{\infty}(\mathbb{D}, X)$. Here $\mathbb{D}= \{z\in \mathbb{C}: |z| <1\}$ denotes the unit disk. We also introduce the following geometric notion of $p$-uniformly $\mathbb{C}$-convexity of order $N$ for a complex Banach space $X$  for some $N \in \mathbb{N}$. 
In this paper, for $p\in [2,\infty)$ and each $N \in \mathbb{N}$, we prove that a complex Banach space $X$ is $p$-uniformly $\mathbb{C}$-convex of order $N$ if, and only if, the $p$-Bohr radius of order $N$ $\widetilde{R}_{p,N}(X)>0$. We also study the $p$-Bohr radius of order $N$ for the Lebesgue spaces $L^q (\mu)$
for  $1\leq p<q<\infty$ or $1\leq q \leq p <2$. Finally, we prove an operator valued analogue of a refined version of Bohr and Rogosinski inequality for bounded holomorphic functions from the unit disk $\mathbb{D}$ into $\mathcal{B(\mathcal{H})}$, where $\mathcal{B(\mathcal{H})}$ denotes the space of all bounded linear operator on a complex Hilbert space $\mathcal{H}$.
\end{abstract}

\maketitle
\pagestyle{myheadings}
\markboth{Vasudevarao Allu, Himadri Halder, and Subhadip Pal}{Bohr and Rogosinski inequalities for operator valued holomorphic functions}
\vspace{-2mm}

\section{Introduction}
Let $H^{\infty}(\mathbb{D}, X)$ be the space of all bounded holomorphic functions from the unit $\mathbb{D}=\left\{z\in \mathbb{C} : |z| <1\right\}$ into a complex Banach space $X$ and we denote $\norm{f}_{H^{\infty}(\mathbb{D}, X)}= \sup_{|z|<1}\norm{f(z)}$. The remarkable discovery of H. Bohr \cite{Bohr-1914} for the functions in $H^{\infty}(\mathbb{D}, \mathbb{C})$ states that 
\begin{thm}
	If $f \in H^{\infty}(\mathbb{D}, \mathbb{C})$ and $f(z)= \sum_{n=0}^{\infty}a_n z^n$ in $\mathbb{D}$, then 
	\begin{align}\label{e-1.1}
		\sum_{n=0}^{\infty}|a_n|r^n \leq \norm{f}_{ H^{\infty}(\mathbb{D}, \mathbb{C})}
	\end{align}
for $|z|=r \leq 1/3$, and the constant $1/3$, referred to as the classical Bohr radius which is sharp.
\end{thm}
The Bohr's theorem has been found a new direction of research when Dixon \cite{Dixon & BLMS & 1995} used it to disprove the conjecture that a non-unital Banach algebra that satisfies the von-Neumann inequality must be isometrically isomorphic to a closed subalgebra of $B(H)$ for some Hilbert space $H$. Since then, Bohr's result got various extensions by several authors. In $1997$, Boas and Khavinson \cite{boas-1997} have established a multidimensional analogue of Bohr inequality in the complete circular domain in $\mathbb{C}^n$. Later other variations of this inequality in the context of several complex variables have been extensively studied by several authors (see \cite{aizn-2000,aizenberg-2001,aizn-2007,Himadri-Vasu-P9,Ayt & Dja & BLMS & 2013,bayart-advance-2014,defant-2006, Liu-Pon-PAMS-2020, paulsen-2002,popescu-2019}). In $2010$, Blasco \cite{Blasco-OTAA-2010} extended the Bohr inequality for Banach spaces and introduced the new notion of Bohr radius for vector valued functions in the unit disk $\mathbb{D}$ (see also \cite{Blasco-Collect-2017}). For more intriguing aspects of Bohr radius for vector valued holomorphic functions in connection with local Banach space theory, we refer to \cite{defant-2003,defant-2008,defant-2012,defant-2018}. Recently, Bhowmik and Das \cite{bhowmik-Arch-Math-2021} have extensively studied a characterization of Banach spaces with nonzero Bohr radius and obtain a necessary and sufficient condition for the occurrence of the Bohr phenomenon for Banach space valued holomorphic functions. For the interesting work on Bohr radius for complex valued holomorphic functions, we refer to \cite{Ahamed-Allu-Halder-P3-2020,Himadri-Vasu-P1,Himadri-Vasu-P2}.
\\[2mm]

For the last two decades,  numerous versions of Bohr inequality has been studied by several authors. In \cite{bene-2004}, B\'{e}n\'{e}teau {\it et al.} have considered the Bohr inequality for the functions in the unit ball of $H^p(\mathbb{D}, \mathbb{C})$ for $p>0$ instead of the functions in $H^{\infty}(\mathbb{D}, \mathbb{C})$. In 2000, Djakov and Ramanujan \cite{Djakov & Ramanujan & J. Anal & 2000} generalized the Bohr inequality \eqref{e-1.1} by introducing the best constant $r_p$ satisfying 
\begin{equation} \label{e-1.2}
	\left(\sum_{n=0}^{\infty} |a_n|^p (r_p)^{np}\right)^{1/p} \leq \norm{f}_{H^{\infty}(\mathbb{D}, \mathbb{C})}
\end{equation}
for any $f(z)= \sum_{n=0}^{\infty}a_n z^n \in H^{\infty}(\mathbb{D}, \mathbb{C})$. Clearly for $p=1$, the radius $r_1 = 1/3$, which is the classical Bohr radius. In view of Haussdorff-Young's inequality, one can easily see that $r_{p}=1$ for $p\geq 2$. 
\vspace{1mm}
Besides Bohr inequality, there is another classical inequality given by W. Rogosinski \cite{rogosinski-1930} which states that, if $f(z)=\sum_{k=0}^{\infty}a_k z^k$ be a holomorphic function on the unit disk $\mathbb{D}$, then for every $N\geq 1$, we have $|\sum_{k=0}^{N}a_k z^k| \leq 1$ in the disk $|z| \leq 1/2$, and the constant $1/2$ is the best possible. 
In recent years, an enormous research has been carried out on refined and improved versions of Bohr inequality (see \cite{Huang-Liu-Ponnu,kayumov-2018-c, Liu-2021}). 
\\[2mm]

In \cite{N-Das-JMAA-2022}, Das proved some new refinements of Bohr inequality as well as Rogosinski inequality. One of the aims of this paper is to study those refined inequalities for operator valued holomorphic functions or, precisely for the holomorphic functions from the unit disk $\mathbb{D}$ into $\mathcal{B(\mathcal{H})}$, where $\mathcal{B(\mathcal{H})}$ is the space of all bounded linear operators on a complex Hilbert space $\mathcal{H}$. Recently, Allu and Halder \cite{Himadri-Vasu-P9} have extended some interesting improved and refined versions of Bohr inequality for bounded holomorphic functions with values in $\mathcal{B(\mathcal{H})}$. 
\vspace{1mm}

In $2019$, Popescu \cite{popescu-2019} proved the following interesting result, which is an  operator valued analogue of the classical theorem of Bohr for  operator valued bounded holomorphic functions in the unit disk $\mathbb{D}$. 
\begin{customthm}{A} \cite{popescu-2019}\label{him-vasu-P9-thm-1.1} 
	Let $f \in H^{\infty}(\mathbb{D},\mathcal{B}(\mathcal{H}))$ be an operator valued bounded holomorphic function with $f(z)=\sum_{n=0}^{\infty}A_{n}z^n$ in $\mathbb{D}$ such that $A_{k} \in \mathcal{B(\mathcal{H})}$ for all $k \in \mathbb{N} \cup \{0\}$ and $A_{0}=a_{0}I$, $a_{0} \in \mathbb{C}$. Then 
	\begin{equation} \label{him-vasu-P9-e-1.1}
		\sum _{n=0}^{\infty} \norm{A_{n}} \, r^n \leq \norm{f}_{H^{\infty}(\mathbb{D},\mathcal{B}(\mathcal{H}))} \,\,\,\,\,\, \mbox{for} \,\,\,\,\, |z|=r \leq \frac{1}{3}
	\end{equation}
	and $1/3$ is the best possible constant. Moreover, the inequality is strict unless $f$ is a constant.
\end{customthm}

Another interesting direction of Bohr phenomenon is getting enriched by considering the Bohr inequality for subordinating families of holomorphic functions in the unit disk $\mathbb{D}$. Indeed, for two holomorphic functions $f$ and $g$ in $\mathbb{D}$, we say that $f$ is subordinate to $g$, denoted by $f \prec g$ , if there exists a holomorphic function $\phi : \mathbb{D} \rightarrow \mathbb{D}$ with $\phi(0)=0$ such that $f(z)= g(\phi(z)),\, z\in \mathbb{D}$. Bhowmik and Das \cite{bhowmik-2018} have established that if $f$ and $g$ are two holomorphic functions in $\mathbb{D}$ with Taylor series expansions $\sum_{k=0}^{\infty}a_n z^n$ and $g(z)=\sum_{k=0}^{\infty}b_n z^n$ respectively such that $f \prec g$, then 
\begin{equation} \label{e-1.3}
	\sum_{k=0}^{\infty}|a_n|r^n \leq \sum_{k=0}^{\infty}|b_n|r^n
\end{equation} 
holds for $|z|=r\leq 1/3$. In \cite{bhowmik-Edinburgh-2021}, Bhowmik and Das extended the inequality \eqref{e-1.3} for holomorphic functions $f:\mathbb{D} \rightarrow \mathcal{B(\mathcal{H})}$.
\vspace{1mm}

We now introduce some basic notations for remaining discussions in our paper. For any $A \in \mathcal{B(\mathcal{H})}$, $\norm{A}$ denotes the operator norm of $A$, and $|A|=(A^*A)^{1/2}$ denotes the absolute value of $A$, where $A^*$ denotes the adjoint of $A$ and $B^{1/2}$ denotes the unique positive square root of a positive operator $B$. Throughout the paper the operators $I$ and $O$ will denote the identity operator and the zero operator respectively. 
\vspace{1mm}

The presentation of the paper is as follows. In Section 2, we discuss the notion of $p$-Bohr radius of order $N$ and its relation with geometry of Banach spaces and prove the key result  Theorem \ref{thm-2.4}. In Section 3, we study the refinements of Bohr and Rogosinski inequality for operator valued holomorphic functions $f: \mathbb{D} \rightarrow \mathcal{B(\mathcal{H})}$.

\vspace{3mm}
\section{Geometric Characterization of $p$-Bohr radius}
In \cite{Blasco-Collect-2017}, Blasco has introduced the $p$-Bohr radius by considering the vector-valued analogue of \eqref{e-1.2}. 
\begin{defn}
	Let $1\leq p<\infty$  and let $X$ be a complex Banach space. We write 
	\begin{align}
		r_p(f, X)= \sup \left\{r\geq 0 : \sum_{n=0}^{\infty} \norm{x_n}^p r^{np} \leq 1\right\},
	\end{align}
where $f(z)= \sum_{n=0}^{\infty}x_n z^n$ with $\norm{f}_{H^{\infty}(\mathbb{D}, X)} \leq 1$ and define the $p$- Bohr radius of $X$ by
\begin{align}
	r_p(X)= \inf \left\{r_p(f, X) : \norm{f}_{H^{\infty}(\mathbb{D}, X)} \leq 1 \right\}.
\end{align} 
\end{defn}
In particular, $r_{p}(\mathbb{C}) \geq p/(p+2)$ (see \cite{Blasco-Collect-2017}). It is easy known that $r^{p_{1}}_{p_{1}}(X) \leq r^{p_{2}}_{p_{2}}(X)$ for any Banach space $X$ and $p_{1} \leq p_{2}$ (see \cite{Blasco-Collect-2017}). The $p$-Bohr radius is very much related to the strong maximum modulus theorem. In \cite{Blasco-Collect-2017}, Blasco has shown that if $r_{p}(X)>0$ for some $1 \leq p < \infty$ then $X$ satisfies the strong maximum modulus theorem for any complex Banach space $X$. The strong maximum modulus theorem has an interesting connection with  certain geometric property of the complex Banach spaces, namely strict $c$-convexity (see \cite{dilworth-1986,globevnik-1975,thorp-PAMS-1967}). Now we introduce the notion of $p$-uniformly $\mathbb{C}$-convexity of order $N$ for a complex Banach space $X$ for some $N \in \mathbb{N}$.
\begin{defn} \label{def-1.2}
	Let $2 \leq p < \infty$. A complex Banach space $X$ is called $p$-uniformly $\mathbb{C}$-convex of order $N$ if there exists a constant $\lambda > 0$ such that 
	\begin{align}\label{e-1.8}
	\left(\norm{x_0}^p + \lambda \norm{x_1}^p + {\lambda}^2 \norm{x_2}^p + \cdots + {\lambda}^N \norm{x_N}^p \right)^{1/p} \leq \max_{\theta \in [0, 2\pi)} \norm{x_0 + \sum_{k=1}^{N}e^{i \theta}x_k}
	\end{align}
	for all $x_0$, $x_1$,$\dots$, $x_N$ $\in X$. We denote $A_{p,N}(X)$, the supremum of such constants $\lambda$ satisfying \eqref{e-1.8}.
\end{defn}
In particular, for $N=1$, Definition \ref{def-1.2} gives the definition of usual $p$-uniformly $\mathbb{C}$-convexity of a complex Banach space $X$. In \cite{Blasco-Collect-2017}, Blasco has proved the following interesting result which gives a characterization of $p$-uniformly $\mathbb{C}$-convexity of a complex Banach space $X$ with the nonzero $p$-Bohr radius. 
\begin{thm} \label{thm-1.2}\cite{Blasco-Collect-2017}
	Let $X$ be a complex Banach space and $p\geq 2$. Then $X$ is $p$-uniformly $\mathbb{C}$-convex if, and only if, the $p$-Bohr radius $r_{p}(X)>0$.
\end{thm}  
In 2003, Blasco and Pavlovic \cite{Blasco-Pavlovic-BLMS-2003} established that the property $p$-uniformly $\mathbb{C}$-convexity of a complex Banach space $X$ can also be expressed in the following way: the existence of a constant $\lambda>0$ such that 
\begin{equation} \label{e-1.6}
\left(\norm{f(0)}^p + \lambda \norm{f'(0)}^p\right)^{1/p} \leq \norm{f}_{H^{\infty}(\mathbb{D},X)}
\end{equation} 
for any $f \in H^{\infty}(\mathbb{D},X)$.
 Motivated by \eqref{e-1.6}, Blasco \cite{Blasco-Collect-2017} has introduced the following notion of $p$-Bohr radius.
\begin{defn} \label{defn-1.3} \cite{Blasco-Collect-2017}
	Let $X$ be a complex Banach space and $1 \leq p <\infty$. We denote 
	\begin{align}\label{e-1.7}
		\widetilde{r}_p(X)= \sup \left\{r>0 : \norm{f(0)}^p + r^p \norm{f'(0)}^p \leq \norm{f}^p_{H^{\infty}(\mathbb{D}, X)}\right\}
	\end{align}
for $f(z)= \sum_{n=0}^{\infty} x_n z^n \in H^{\infty}(\mathbb{D}, X)$. 
\end{defn}
The Definition \ref{defn-1.3} leads us to consider the following notion of  $p$-Bohr radius of order $N$.   
Further, we show that this notion characterizes the property $p$-uniformly $\mathbb{C}$-convexity of a complex Banach space $X$. 
\begin{defn}\label{defn-1.4}
	Let $X$ be a complex Banach space and $1 \leq p < \infty$. For each $N \in \mathbb{N}$, we denote
	\begin{align}
		\widetilde{R}_{p,N}(f, X) = \sup \left\{r \geq 0 : \norm{x_0}^p + \sum_{k=1}^{N}\norm{x_k}^p r^{pk} \leq 1\right\},
	\end{align}
where $f(z)= \sum_{k=0}^{\infty}x_k z^k$ with $\norm{f}_{H^{\infty}(\mathbb{D}, X)} \leq 1$.
We define the $p$-Bohr radius of order $N$ by
\begin{align*}
	\widetilde{R}_{p,N}(X) &= \inf \left\{\widetilde{R}_{p,N}(f, X) : \norm{f}_{H^{\infty}(\mathbb{D}, X)} \leq 1 \right\}\\
	& = \sup \left\{r\geq 0: \left(\sum_{k=0}^{N}\norm{x_k}^p r^{pk}\right)^{1/p} \leq \norm{f}_{H^{\infty}(\mathbb{D}, X)}\right\}.
\end{align*}
\end{defn}
\noindent In particular, for $N=1$, we see that $\widetilde{R}_{p,1}(X)=\widetilde{r}_{p}(X)$. To compute the exact value of $\widetilde{R}_{p,N}(X)$ for any complex Banach space $X$ is a very difficult task, in general. In \cite{Blasco-Collect-2017} Blasco has achieved the precise value as $\widetilde{R}_{p,1}(\mathbb{C})= \inf_{0<a<1} {(1-a^p)}^{1/p}/{(1-a^2)}$.
\begin{rem} \label{rem-1.1}
	Let $1 \leq p < \infty$. Then 
	\begin{equation}
		r_p(X) \leq \widetilde{R}_{p,N}(X) \leq \widetilde{r}_p(X).
	\end{equation}
Indeed, since 
$\norm{x_0}^p + r^p \norm{x_1}^p \leq \norm{x_0}^p + \sum_{k=1}^{N}\norm{x_k}^p r^{pk} \leq \norm{x_0}^p + \sum_{k=1}^{\infty} \norm{x_k}^p r^{pk}.$
\end{rem}
We first obtain a lower estimate of $\widetilde{R}_{p,N}(\mathbb{C})$ for any $N \in \mathbb{N}$.
\begin{thm} \label{thm-2.1}
Let $1 \leq p < \infty$ and $N \in \mathbb{N}$. Then $\widetilde{R}_{p,N}(\mathbb{C}) \geq r^{*}_{N}$, where $r^{*}_{N}$ is the positive root of the equation 
\begin{equation}
r^p + r^{2p} + \cdots + r^{Np} - \xi _{p}=0,
\end{equation}
provided $\xi_{p} < N$. Here $\xi_{p}= \inf _{0<a<1} {(1-a^p)^{1/p}/(1-a^2)}$.
\end{thm}
\begin{pf}  
	Let $f(z)= \sum_{k=0}^{\infty} a_{k} z^k \in H^{\infty}(\mathbb{D},\mathbb{C})$ with $\norm{f}_{H^{\infty}(\mathbb{D},\mathbb{C})} \leq 1$. From Wiener's lemma, we have $|a_{n}| \leq 1- |a_{0}|^2$ for all $n \geq 1$. Therefore, we have 
	\begin{align}
		|a_{0}|^p + \sum_{k=1}^{N} |a_{k}|^p r^{pk} & \leq |a_{0}|^p + (1-|a_{0}|^2)^p\sum_{k=1}^{N} r^{pk} 
		\leq |a_{0}|^p + (1-|a_{0}|^2)^p \, r^p  \frac{1-r^{pN}}{1-r^p},
	\end{align} 
	which is less than or equals to $1$ if 
	\begin{equation}\label{e-2.10}
		r^p \, \frac{1-r^{pN}}{1-r^p} \leq \frac{1-|a_{0}|^p}{(1-|a_{0}|^2)^p}.
	\end{equation}
	Let $\Psi_{p,N}(r)= r^p+r^{2p}+\cdots+r^{Np}-\xi_{p}$ . 
	Then the inequality \eqref{e-2.10} holds good for $r\leq r_N^{*}$, where $r_N^{*}$ is a positive root of the equation $\Psi_{p,N}(r)=0$ in $(0,1)$.
	We now show the uniqueness of the root $r_N^{*}$.	
	Here 
	\begin{align*}
		\Psi_{p,N}^{'}(r)= \sum_{k=1}^N pk r^{pk-1}.
	\end{align*}
	Clearly $\Psi_{p,N}^{'}(r)>0$ for $r\in (0,1)$. Therefore $\Psi_{p,N}^{'}(r)$ is an increasing function. Further, $\Psi_{p,N}(0)=-\xi_p<0$ and $\Psi_{p,N}(1)=N-\xi_{p}>0$.
	Hence the function $\Psi_{p,N}$ has the unique root in $(0,1)$.
	This completes the proof.
\end{pf}

Next, we establish a result on both bounds of the $p$-Bohr radius in terms of $p$-Bohr radius of order $N$, which generalizes \cite[Theorem 2.4]{Blasco-Collect-2017}. The following inequality is useful to prove one of the main results of our paper.
\begin{thm} \label{thm-2.2}
	Let $p\geq1$ and $N \in \mathbb{N}$. Then 
	\begin{equation}\label{e-2.2}
		\frac{\widetilde{R}_{p,N}(X)}{\left(1+\widetilde{R}_{p,N}(X)^p\right)^{1/p}} \leq r_p(X) \leq \widetilde{R}_{p,N}(X).
	\end{equation}
\end{thm}
\begin{pf} 
	Let $f\in H^{\infty}(\mathbb{D},X) $ with norm $1$ and $f(z)= \sum_{n=0}^{\infty} x_n z^n$.
	Let $n\in \mathbb{N}$ and consider $\xi= e^{{2\pi i}/{n}}$ and define $g(z)= \frac{1}{n} \sum_{j=1}^{n}f(\xi^j z)$ (see \cite{Bohr-1914}). Using the fact that $\sum_{j=1}^{n}\xi^j=0$ we obtain 
	\begin{align}\nonumber
		g(z)= x_0 + x_nz^n+ x_{2n}z^{2n}+ \cdots \in H^{\infty}(\mathbb{D},X).
	\end{align}
	Clearly $\norm{g}_{H^{\infty}(\mathbb{D},X)} \leq 1$, indeed 
	\begin{equation}\nonumber
		\norm{g}_{H^{\infty}(\mathbb{D},X)}\leq \frac{1}{n} \sum_{j=1}^{n}\norm{f(\xi_jz)}_X=1.
	\end{equation}
	Then by the definition of $\widetilde{r}_p(X)$, we have 
	\begin{equation*}
		\norm{g(0)}^p + \widetilde{r}_p(X)^p \norm{g^{'}(0)} \leq 1, 
	\end{equation*}
	which gives
	\begin{equation*}
		\norm{g'(0)}^p \leq \widetilde{r}_p(X)^{-p}(1- \norm{g(0)}^p).
	\end{equation*}
	Thus we obtain the following estimate
	\begin{align}\label{e-2.12}
		\norm{x_n}^p & \leq \tilde{r}_p(X)^{-p}(1- \norm{x_0}^p) \\ \nonumber & \leq \widetilde{R}_{p,N}(X)^{-p}(1- \norm{x_0}^p), \,\,\, n\geq 1. 
	\end{align}	
	Therefore using \eqref{e-2.12} we obtain
	\begin{align} \nonumber
		\norm{x_0}^p + \sum_{n=1}^{\infty} \norm{x_n}^pr^{np} & \leq \norm{x_0}^p + \widetilde{R}_{p,N}(X)^{-p}(1- \norm{x_0}^p)\left(\sum_{n=1}^{\infty}r^{np}\right) \\ \nonumber & \leq \norm{x_0}^p + \widetilde{R}_{p,N}(X)^{-p}(1- \norm{x_0}^p) \frac{r^p}{1-r^p} \\ \nonumber & \leq \max \left\{1, \widetilde{R}_{p,N}(X)^{-p}\frac{r^p}{1- r^p}\right\}.
	\end{align}
	Now we choose $r$ in such a way that $\widetilde{R}_{p,N}(X)^{p}= \frac{r^p}{1-r^p}$. Thus we obtain 
	\begin{align}\nonumber
		\frac{\widetilde{R}_{p,N}(X)}{\left(1+\widetilde{R}_{p,N}(X)^p\right)^{1/p}} \leq r_p(X).
	\end{align} 
	In view of Remark \ref{rem-1.1}, we obtain $r_p(X) \leq \widetilde{R}_{p,N}(X) $. This completes the proof.
\end{pf}

From \eqref{e-2.2} in Theorem \ref{thm-2.2} we obtain the following lower estimate of $r_p(\mathbb{C})$.
\begin{cor}
	Let $1<p<2$. Then 
	\begin{align}
		r_p(\mathbb{C}) \geq \left\{1+ \widetilde{R}_{p,N}(\mathbb{C})^{-p}\right\}^{-1/p}.
	\end{align}
\end{cor}
It is known that $p$-uniformly $\mathbb{C}$-convexity of a complex Banach space is strongly related to $p$-Bohr radius. Now we state  the following result where we obtain both lower and upper bounds of $\widetilde{R}_{p,N}(X)$ in terms of $A_{p,N}(X)$, which has a significant importance on our next result.
\begin{thm} \label{thm-2.3}
	Let $p \geq 2$ and $N \in \mathbb{N}$. Then 
	\begin{align}\label{e-2.4}
		\frac{A_{p,N}^{\frac{1}{p}}(X)}{2N} \leq \widetilde{R}_{p,N}(X) \leq A_{p,N}^{\frac{1}{p}}(X).
	\end{align}
\end{thm}
\begin{pf}  
	Let $f\in H^{\infty}(\mathbb{D}, X)$ with $\norm{f}_{H^{\infty}(\mathbb{D}, X)} \leq 1$ and $f(z)= \sum_{n=0}^{\infty} x_n z^n$. Let $y= e^{{2\pi i}/n}$ for $n \in \mathbb{N}$ and we define $h(z)= 1/n\sum_{k=1}^{n}f(y^k z)$.
	Since $\sum_{k=1}^{n}y^k=0$, we obtain a function $h\in H^{\infty}(\mathbb{D}, X)$, where $h(z)= \sum_{k=0}^{\infty} x_{nk} z^{nk}$ with $\norm{h}_{H^{\infty}(\mathbb{D}, X)}\leq 1$.
	Let $\xi \in X^{*}$ with norm $1$, where $X^*$ denotes the dual space of $X$. Then the composition function $ \xi \circ h= \innpdct{\xi, h} \in H^{\infty}(\mathbb{D}, \mathbb{C}) $ with $\norm{\innpdct{\xi, h}} \leq \norm{\xi}\norm{h}\leq 1$. In view of the Schwarz-Pick lemma, we obtain 
	\begin{align}\nonumber
		|{\innpdct{\xi, h'(0)}}| & \leq 1- |\innpdct{\xi, h(0)}|^2 \\ \nonumber & \leq 2\left(1-|\innpdct{\xi, h(0)}|\right).
	\end{align}
	Therefore we have
	\begin{align}\nonumber
		|\innpdct{\xi, x_n}| \leq 2\left(1- |\innpdct{\xi, x_0}|\right), \quad n\geq 1.
	\end{align}	
	Hence for any $N \in \mathbb{N}$
	\begin{align}\nonumber
		\sum_{k=1}^{N}|\innpdct{\xi, x_k}| \leq 2N\left(1- |\innpdct{\xi,x_0}|\right),
	\end{align}
	that is,
	\begin{align} \label{e-3.7-a}
		|\innpdct{\xi, x_0}| + \frac{1}{2N}\sum_{k=1}^{N}|\innpdct{\xi,x_k}| \leq 1.
	\end{align}
	Consequently for any $\theta \in [0, 2\pi)$,
	\begin{align}\nonumber
		\norm{x_0 + \frac{e^{i\theta}}{2N}\sum_{k=1}^{N} x_k} &=  \sup_ {\norm{\xi}=1}\left|\innpdct{\xi, x_0} + \frac{e^{i\theta}}{2N} \sum_{k=1}^{N}\innpdct{\xi, x_k}\right| \\ \nonumber & \leq \sup_{\norm{\xi}=1} \left\{| \innpdct{\xi, x_0}| + \frac{1}{2N}\sum_{k=1}^{N}|\innpdct{\xi, x_k}|\right\} \\ \nonumber & \leq \sup_{\norm{\xi}=1} \left\{|\innpdct{\xi, x_0}|+ \frac{1}{2N}2N\left(1-|\innpdct{\xi, x_0}|\right)\right\} \quad (\mbox{using} \quad \eqref{e-3.7-a})\\ \nonumber & =1.
	\end{align}
	Hence we have
	\begin{align}\nonumber
		\norm{x_0}^p + \frac{A_{p,N}(X)}{(2N)^p}\norm{x_1}^p+ \frac{A_{p,N}^2(X)}{(2N)^p}\norm{x_2}^p+ \cdots + \frac{A_{p,N}^N(X)}{(2N)^p}\norm{x_N}^p \leq 1,
	\end{align}
	which shows that ${A_{p,N}^{1/p}(X)}/{2N} \leq \widetilde{R}_{p,N}(X)$. To prove the other inequality, let $f \in H^{\infty}(\mathbb{D}, X)$, then we have
	\begin{align}\label{e-3.7}
		\left(\sum_{k=0}^{N}\norm{x_k}^p r^{pk}\right)^{1/p} \leq \norm{f}_{H^{\infty}(\mathbb{D}, X)}.
	\end{align}
	In view of the definition of $A_{p,N}(X)$, we have $r \leq A^{1/p}_{p,N}(X)$ and hence by taking supremum over all such $r$ satisfying \eqref{e-3.7} we obtain $\widetilde{R}_{p,N}(X) \leq A^{1/p}_{p,N}(X)$. This completes the proof.
\end{pf}

We are now ready to state one of the main results of our paper which characterizes the notion of $p$-uniformly $\mathbb{C}$-convexity property of a complex Banach space with the nonzero $p$-Bohr radius of order $N$. 

\begin{thm} \label{thm-2.4} 
	Let $X$ be a complex Banach space and $p\geq 2$. Then $X$ is $p$-uniformly $\mathbb{C}$-convex of order $N$ if, and only if, $\widetilde{R}_{p,N}(X)>0$.
\end{thm}
\begin{pf}  
	In view of \eqref{e-2.2} and \eqref{e-2.4}, we define the function $h_p(t)$ by
	$$h_p(t)= \frac{t}{(1+t^p)^{1/p}}, \quad t\in [a,b], \quad \mbox{where}\quad	a=\frac{A_{p,N}^{1/p}(X)}{2N} \quad \mbox{and}\quad b=A_{p,N}^{1/p}(X).$$
	A simple computation shows that $$h_p'(t)= \frac{1}{(1+t^p)^{\frac{1}{p}+1}}>0, \quad  t \in [a, b].$$\\
	Thus $h_p$ is an increasing function in $[a, b]$, and hence
	it attains its minimum at $a$. Therefore, we have 
	\begin{align}\nonumber
		h_p\left(\frac{A_{p,N}^{1/p}(X)}{2N}\right) \leq h_p\left(\widetilde{R}_{p,N}(X)\right).  
	\end{align}
	Equivalently,
	\begin{align}\label{e-3.8}
		\frac{A_{p,N}^{1/p}(X)}{\left(A_{p,N}(X)+ (2N)^p\right)^{1/p}} \leq \frac{\widetilde{R}_{p,N}(X)}{\left(1+\widetilde{R}_{p,N}^p(X)\right)^{1/p}}.
	\end{align}
	Combining \eqref{e-2.2}, \eqref{e-2.4} and \eqref{e-3.8} we obtain 
	\begin{align}\label{e-3.9}
		\frac{A_{p,N}^{1/p}(X)}{\left(A_{p,N}(X)+ (2N)^p\right)^{1/p}} \leq r_p(X) \leq A_{p,N}^{1/p}(X).
	\end{align}
	Now the $p$-uniformly $\mathbb{C}$-convexity of order $N$ of a complex Banach space $X$ is equivalent to the condition that $A_{p,N}(X)>0$. Therefore, we obatin $r_p(X)>0$ and  from Remark \ref{rem-1.1} it follows that $\widetilde{R}_{p,N}(X)>0$. Conversely, we assume $\widetilde{R}_{p,N}(X)>0$. Then the inequality \eqref{e-2.4} gives $A_{p,N}(X)>0$. Consequently, $X$ is $p$-uniformly $\mathbb{C}$-convex of order $N$. This completes the proof.
\end{pf}

Combining \eqref{e-2.2} and \eqref{e-2.4} we obtain the following lower estimate of $r_2(X),$ where $X$ is $2$-uniformly $\mathbb{C}$-convex of order $N$.
\begin{cor}
	If $X$ is $2$-uniformly $\mathbb{C}$-convex of order $N$ then $r_2(X) \geq \sqrt{\frac{A_{p,N}(X)}{A_{p,N}(X)+ 4N^2}}$.
\end{cor}
Now we study the $p$-Bohr radius of order $N$ in some infinte dimensional spaces. Motivated by the work of Blasco \cite{Blasco-Collect-2017}, we now establish the $p$-Bohr radius of order $N$ for the Lebesgue spaces $L^q(\mu)$.\\
\vspace{2mm}

Let $X=L^q(\mu)$, where $\mu$ is a measure and $1 \leq q\leq \infty$. In this case we see that the same situation arises for $X$ whenever $(\Omega, \Sigma, \mu)$ is a measure space and there exists a sequence of pairwise disjoint measurable sets $\{A_n\}$ with $0< \mu(A_n)<\infty$. In fact, Blasco \cite{Blasco-Collect-2017} has proved the following result. 
\begin{thm}\cite{Blasco-Collect-2017}\label{thm-2.5}
	Let $(\Omega, \Sigma, \mu)$ be a measure space such that there exists a couple of disjoint measurable sets $A, B \in \Sigma$ with $0 < \mu(A), \mu(B) < \infty$. Then 
	\begin{equation}\label{e-2.5}
		\widetilde{r}_p(L^q(\mu)) = 0, \quad 1 \leq p<q< \infty \quad or \quad  1\leq q \leq p<2.
	\end{equation}
\end{thm}
Now we obtain the following analogue for $p$-Bohr radius of order $N$.
\begin{thm} \label{thm-2.6}
	Let $1 \leq p, q < \infty$ and let $(\Omega, \Sigma, \mu)$ be a measure space such that there exists a couple of disjoint measurable sets $A, B \in \Sigma$ with $0< \mu(A), \mu(B)< \infty$. Then for $N \geq 1$,
	\begin{equation*}
		\widetilde{R}_{p,N}(L^q(\mu))=0, \quad 1 \leq p<q< \infty \quad or \quad  1\leq q \leq p<2.
	\end{equation*}
\end{thm}
\begin{pf} 
	From Remark \ref{rem-1.1}, we have the inequality $\widetilde{R}_{p,N}(X) \leq \widetilde{r}_p(X)=\widetilde{R}_{p,1}(X)$.
	Therefore, from \eqref{e-2.5} we obtain
	\begin{equation*}
		\widetilde{R}_{p,N}(L^q(\mu))\leq \widetilde{r}_p(L^q(\mu))=0
	\end{equation*}
	for $1\leq p<q<\infty$ or $1\leq q \leq p <2$ and for $N \geq 1$.
	In view of the Definition \ref{defn-1.4} we have, $\widetilde{R}_{p,N}(L^q(\mu)) \geq 0$. Consequently, $\widetilde{R}_{p,N}(L^q(\mu))=0$ for $1\leq p<q<\infty$ or $1\leq q \leq p <2$ and for $N \geq 1$. This completes the proof.
\end{pf}



\vspace{-5mm}

 \vspace{5mm}
\section{Operator valued analogue of refinements of the Bohr and Rogosinski inequality}
In this section we study an operator valued analogue of refined form of the classical Bohr and Rogosinski inequalities. 
For any holomorphic function $f: \mathbb{D} \rightarrow \mathcal{B}(\mathcal{H})$ with $f(z) = \sum_{n=0}^{\infty} A_n z^n$ in $\mathbb{D}$ and $A_n \in \mathcal{B}(\mathcal{H})$ for $n \in \mathbb{N} \cup \{0\}$, we define the majorant series $M_r(f)$ for $f$ by $$ M_r(f)= \sum_{n=0}^{\infty} \norm{A_n}|z|^n  \quad \mbox{for} \,\, |z|=r \in [0,1).$$ 
 We now discuss some properties of $M_r(f)$. Let $f,g : \mathbb{D} \rightarrow \mathcal{B}(\mathcal{H})$ be two holomorphic functions with power series expansions $f(z)=\sum_{n=0}^{\infty} A_n z^n$ and $g(z)=\sum_{n=0}^{\infty} B_n z^n$ respectively in $\mathbb{D}$, where $A_n, B_n \in \mathcal{B}(\mathcal{H}) $ and $n \in \mathbb{N} \cup \{0\}$. Then we have $M_r(f + g) \leq M_r(f) + M_r(g)$. Indeed, 
 
 \begin{equation}\label{e-4.1}
 	M_r(f+g) = \sum_{n=0}^{\infty} \norm{A_n + B_n} r^n \leq \sum_{n=0}^{\infty} \norm{A_n}r^n + \sum_{n=0}^{\infty} \norm{B_n}r^n = M_r(f) + M_r(g).
 \end{equation}

 Moreover, using the above inequality \eqref{e-4.1} we can prove that if $F(z) = \sum_{j \in \mathbb{Z}} f_{j}(z)$ is holomorphic in $\mathbb{D}$, where each $f_{j}$ is holomorphic in $\mathbb{D}$ for each $j \in \mathbb{Z}$, then $M_r(F)\leq \sum_{j \in \mathbb{Z}} M_r(f_{j})$. Further, we observe that $M_r(\alpha f) = |\alpha | M_r(f)$ for any $\alpha \in \mathbb{C}$, and $M_r(z^k f)= r^k M_r(f)$. In the case of holomorphic function $(fg)(z)= \sum_{n=0}^{\infty} A_n(z^n g(z))$ we have
 \begin{equation}
 	M_r(fg)\leq \sum_{n=0}^{\infty} \norm{A_n}r^n M_r(g)= M_r(f)M_r(g).
 \end{equation}
In 2021, Bhowmik and Das \cite{bhowmik-Edinburgh-2021} proved the following lemma for operator valued analogue of subordination for majorant series.

\begin{lem} \cite{bhowmik-Edinburgh-2021}
	Let $f, g : \mathbb{D} \rightarrow \mathcal{B}(\mathcal{H})$ be holomorphic functions with expansions $f(z)= \sum_{n=0}^{\infty}A_n z^n$ and $g(z)=\sum_{n=0}^{\infty}B_n z^n$ respectively such that $f \prec g$. Then 
	\begin{equation}
		\sum_{n=0}^{\infty}\norm{A_n}r^n \leq \sum_{n=0}^{\infty}\norm{B_n}r^n \quad for \,\, |z|=r \leq \frac{1}{3}.
	\end{equation}
\end{lem}
We now obtain the following operator valued analogue of a refined Bohr inequality which has been proved recently by Das \cite{N-Das-JMAA-2022}. 
\begin{thm}
	Let $f: \mathbb{D} \rightarrow \mathcal{B(\mathcal{H})}$ be a non-constant bounded holomorphic function with the expansion $f(z)= \sum_{n=0}^{\infty} A_n z^n$ in $\mathbb{D}$ such that $A_k \in \mathcal{B(\mathcal{H})}$ for all $k \in \mathbb{N} \cup \{0\}$ and $A_{0} = \alpha I$, $\alpha \in \mathbb{C}$ with $|\alpha| \leq 1$. If $\norm{f(z)}\leq 1$ in $\mathbb{D}$, then we have 
	\begin{equation}\label{e-5.4}
		\sum_{n=0}^{\infty} \norm{A_n} r^n + (1- |\alpha|) \frac{1-r(1+2|\alpha|)}{1-r|\alpha|} + \frac{r(1-|\alpha|^2)(1-G)}{(1-r|\alpha|)(1-r|\alpha|G)} \leq 1
	\end{equation}
for $|z|=r\leq 1/3$, where
\begin{equation*}
	G= \beta + \frac{r(1-{\beta}^2)}{1-r\beta}, \,\,\,\, \beta = \frac{\norm{A_1}}{1-|\alpha|^2}.
\end{equation*}
\end{thm}
\begin{pf}
	Let $\psi(z)=\sum_{k=1}^{\infty}b_k z^k : \mathbb{D} \rightarrow \mathbb{D}$ be a holomorphic function. By the well-known result of subordination we can write $$\psi(z)/z \prec \frac{b_1 - z}{1- \overline{b_1}z},\,\, z \in \mathbb{D}.$$ By applying \cite[Lemma 1]{bhowmik-2018}, we obtain 
	\begin{equation*}
		\sum_{n=1}^{\infty} |b_n|r^n \leq r \bigg(|b_1| + \frac{r(1-|b_1|^2)}{1-r|b_1|}\bigg) \,\,\, \mbox{for} \,\,\, |z|=r\leq \frac{1}{3}.
	\end{equation*}
	The assumption $\norm{f(z)}\leq 1$ guarantees the existence of a biholomorphic function $g: \mathbb{D} \rightarrow \mathcal{B(\mathcal{H})}$ with $\norm{g(z)}\leq 1 $ such that $f \prec g$. Without loss of generality, let us take the biholomorphic function 
	\begin{equation*}
		g(z)= \bigg(\frac{\alpha - z}{1- \overline{\alpha}z}\bigg) I, \,\, z\in \mathbb{D}.
	\end{equation*}
	Therefore, we have
	\begin{equation}\label{e-5.5}
		f(z)= \bigg(\frac{\alpha - \psi(z)}{1- \overline{\alpha}\psi(z)}\bigg) I = \alpha I - (1-|\alpha|^2)\sum_{n=1}^{\infty}(\overline{\alpha})^{n-1}\psi^{n}(z) I, \,\, z\in \mathbb{D}
	\end{equation}
for some $\psi$ discussed above. From \eqref{e-5.5}, it is evident that $|b_1|=\norm{A_1}/(1-|\alpha|^2)$ and $\beta= |b_1|$.\\

Now based on the discussions on properties of majorant series made at the beginning of this section and from \eqref{e-5.5} we obtain
\begin{align*}
	\sum_{n=0}^{\infty} \norm{A_n}r^n 
	&= \norm{\alpha}I + M_r\bigg(\sum_{n=1}^{\infty}(1-|\alpha|^2)(\overline{\alpha})^{n-1}\psi^n\bigg) \\ 
	& \leq |\alpha| + (1-|\alpha|^2) \sum_{n=1}^{\infty}|\alpha|^{n-1}\big(M_r(\psi^n)\big) \\
	 & \leq  |\alpha| + (1-|\alpha|^2) \sum_{n=1}^{\infty}|\alpha|^{n-1}\bigg(\sum_{k=1}^{\infty}|b_k|r^k\bigg)^n \\ 
	 & \leq  |\alpha| + (1-|\alpha|^2) \sum_{n=1}^{\infty}|\alpha|^{n-1}\bigg(r\big(\beta + \frac{r(1-\beta^2)}{1-r\beta}\big)\bigg)^n \\ 
	 & = |\alpha| + \frac{r(1-|\alpha|^2)}{1-r|\alpha|}- \frac{r(1-|\alpha|^2)(1-G)}{(1-r|\alpha|)(1-r|\alpha|G)} \\ &:= 1- \Psi(r), 
\end{align*}
where $$\Psi(r)= (1-|\alpha|) \big(\frac{1-r(1+2|\alpha|)}{1-r|\alpha|}\big)+\frac{r(1-|\alpha|^2)(1-G)}{(1-r|\alpha|)(1-r|\alpha|G)} \,\,\,\mbox{for}\,\, r \leq 1/3.$$

It is evident that $G\leq 1$ for $r\leq 1/3$. Therefore, we have $\Psi(r) \geq 0$ and hence, the inequality \eqref{e-5.4} holds for $r \leq 1/3$. It is to be noted that if the inequality \eqref{e-5.4} does not give an improved version of the inequality \eqref{him-vasu-P9-e-1.1} for some $r \neq 0$, then
$$\Psi(r)=0 \,\,\,\,\mbox{implies}\,\, (1-|\alpha|) \big(\frac{1-r(1+2|\alpha|)}{1-r|\alpha|}\big)=\frac{r(1-|\alpha|^2)(1-G)}{(1-r|\alpha|)(1-r|\alpha|G)}=0.$$
The first equation gives $|\alpha|=1$ or, $|\alpha|=(1/2)(1/r-1)$ for $r\leq 1/3$
and later one gives $|\alpha| =1$. So, both cases imply $f$ is a constant function, which contradicts our assumption. This completes the proof.
\end{pf}

Now we want to study an operator valued analogue of an refinement of Rogosinski inequality which has been proved recently by Das \cite{N-Das-JMAA-2022}. Let $S(x)$ be the infinte series 
$$S(x)= \bigg(\frac{x}{2}+ \sum_{k=2}^{\infty}\frac{1.3.\cdots (2k-3)}{2^k k!}x^k\bigg)=1-\sqrt{1-x}.$$
Further assume $x_k$ and $y_k$ are real numbers with $|x_k| + |y_k| \neq 0$ for all $1 \leq k \leq n$. Then
$$\left(\sum_{k=1}^{n}x_k y_k\right)^2 \leq \sum_{k=1}^{n}(x_k^2 + y_k^2)\sum_{k=1}^{n} \frac{x_k^2 y_k^2}{x_k^2+ y_k^2} \leq \sum_{k=1}^{n}x_k^2 \, \sum_{k=1}^{n}y_k^2.\quad (Milne's \,\, inequality)$$
We shall prove the following theorem for the holomorphic functions $f \in H^{\infty}(\mathbb{D}, \mathcal{B(\mathcal{H})})$, where $\mathcal{B(\mathcal{H})}$ is a complex Hilbert space.
\begin{thm}
	Let $f: \mathbb{D} \rightarrow \mathcal{B(\mathcal{H})}$ be a non-constant holomorphic function such that $\norm{f(z)}\leq 1$, where $f(z)=\sum_{n=0}^{\infty}A_n z^n$. Then the following inequalities  hold for $|z|=r \leq 1/2$ :
	\item[(a)] $$\norm{\sum_{n=0}^{N}A_n z^n} + (1-2r) r^{N-1}N S(p(N-1))+ (1-r)^2 \sum_{n=0}^{N-2}r^n(n+1) S(p(n))  + r^N (N+1)S(p(N)) \leq 1$$
	for all $N \in \mathbb{N}$, where
	$$p(n)=\frac{1}{n+1}\sum_{l=1}^{\infty}\norm{\sum_{m=0}^{n}A_{m+l}e^{imt}}^2,\,\, n \geq 0\,\, \mbox{and}\,\, t\in [0, 2\pi) . $$
	
	\item[(b)] $$\norm{\sum_{n=0}^{N}A_n z^n} + (1-2r)r^{N-1}NS(q(N-1))+ (1-r)^2 \sum_{n=0}^{N-2} r^n(n+1)S(q(n)) + r^N(N+1)S(q(N)) \leq 1 $$
	for all $N \in \mathbb{N}$, where 
	$$q(n)= \frac{1}{n+1}\sum_{l=0}^{n}\bigg(\frac{1-\norm{R_{l}}^2}{1+ \norm{R_l}^2}\bigg),\,\, n\geq 0 \,\, \mbox{and }\, R_l= \sum_{u=0}^{l}A_u e^{iut}, \, t\in [0, 2\pi).$$
\end{thm}
\begin{pf} (a)
Let $f(z)=\sum_{n=0}^{\infty}A_n z^n$. Then we have $f(ze^{it})=\sum_{n=0}^{\infty}A_n e^{int}z^n$ for any $t\in [0, 2\pi)$ and $z \in \mathbb{D}$. Now for any $n \geq 1$, we obtain
\begin{align*}
	f(ze^{it})\bigg(\sum_{l=1}^{n}z^{l}\bigg) &=  \bigg(\sum_{m=0}^{\infty}A_m e^{imt}z^m\bigg)\bigg(\sum_{l=1}^{n}z^{l}\bigg) \\ & = \sum_{l=0}^{n}\bigg(\sum_{m=0}^{l}A_m e^{imt}\bigg)z^l + \sum_{l=n+1}^{\infty}\bigg(\sum_{m=l-n}^{l}A_me^{imt}\bigg)z^l \\ & = \sum_{l=0}^{n} \bigg(\sum_{m=0}^{l}A_m e^{imt}\bigg)z^l + \sum_{l=1}^{\infty} \bigg(\sum_{m=0}^{n}A_{m+l}e^{i(m+l)t}\bigg)z^{l+n}, \,\, z \in \mathbb{D}.
\end{align*}
Therefore, for any $z=re^{i\xi}$, $\xi \in [0, 2\pi)$ we have the following relation:
\begin{align*}
	\sum_{l=0}^{n}\norm{R_l}^2 r^{2l} + \sum_{l=1}^{\infty}\norm{T_l}^2 r^{2(l+n)} & = \int_{\xi=0}^{2 \pi} \norm{f(ze^{it})(1+z+\cdots+z^n)}^2\, \frac{d\xi}{2\pi} \\ & \leq \int_{0}^{2\pi} |(1+z+ \cdots + z^n)|^2 \, \frac{d\xi}{2\pi} \\ & \leq \sum_{l=0}^{n} r^{2l},
\end{align*}
where $R_l= \sum_{m=0}^{l}A_m e^{imt}$ and, $T_l=\sum_{m=0}^{n}A_{m+l}e^{i(m+l)t}$.
Letting $r \rightarrow 1^{-}$, we obtain 
\begin{equation}\label{e-5.6}
	\sum_{l=0}^{n}\norm{R_l}^2 + \sum_{l=1}^{\infty} \norm{T_l}^2 \leq n+1.
\end{equation}
Let us now define the operator $H_n := \sum_{l=0}^{n}R_l$ for $n \geq 0$ and, $H_n=O$, the zero operator for $n <0$ in $\mathcal{B(\mathcal{H})}$. Using the inquality \eqref{e-5.6} and the Cauchy-Schwartz inequality we obtain
\begin{equation*}
	\norm{H_n} \leq \sum_{l=0}^{n} \norm{R_l} \leq \sqrt{n+1} \sqrt{n+1 - \sum_{l=1}^{\infty}\norm{T_l}^2} = (n+1)(1- S(p(n))),
\end{equation*}
where 
\begin{align*}
	p(n) & = \frac{1}{n+1} \sum_{l=1}^{\infty}\norm{T_l} \\ & = \frac{1}{n+1 }\sum_{l=1}^{\infty} \norm{\sum_{m=0}^{n}A_{m+l}e^{i(m+l)t}}=\frac{1}{n+1 }\sum_{l=1}^{\infty} \norm{\sum_{m=0}^{n}A_{m+l}e^{imt}}
\end{align*}
exactly as defined in part $(a)$ of the statement of this theorem. From \eqref{e-5.6} it is clear that $0 \leq p(n) \leq 1$ for all $n \geq 1$. Hence the series $p(n)$ converges. Now for any $N \in \mathbb{N}$ and $z=re^{it} \in \mathbb{D}$ where $r\leq 1/2$, we obtain
\begin{align*}
	\norm{\sum_{n=0}^{N}A_n z^n} &= \norm{\sum_{n=0}^{N-2}(H_n - 2H_{n-1} + H_{n-2}} \\ & = \norm{(1-r)^2 \sum_{n=0}^{N-2}r^nH_n + (1-2r)r^{N-1}H_{N-1} + r^N H_N} \\ & \leq (1-r)^2 \sum_{n=0}^{N-2}r^n \norm{H_n} + (1-2r)r^{N-1}\norm{H_{N-1}} + r^N \norm{H_N} \\ & \leq (1-r)^2 \sum_{n=0}^{N-2} r^n (n+1)(1-S(p(n))) \\ & \quad \quad + r^{N-1}(1-2r)N(1-S(p(N-1))) \\ & \quad \quad + r^N(N+1)(1-S(p(N))) = 1- \Gamma (r),
\end{align*}
where 
$$ 	\Gamma (r) = (1-2r)\sum_{n=0}^{N-1}r^n (n+1) S(p(n)) + \sum_{n=2}^{N}r^n(n-1)S(n-2) + r^N(N+1)S(p(N)) .$$
Therefore, the inequality in part (a) follows. \\ [2mm]
(b) From \eqref{e-5.6} we have $\sum_{l=0}^{n}\norm{R_l}^2 \leq n+1$, along with the Cauchy-Schwartz inequality yields $\sum_{l=0}^{n}\norm{R_l}\leq n+1$. Now using the Milne's inequality we obtain
\begin{align*}
	\bigg(\sum_{l=0}^{n}\norm{R_l}\bigg)^2 & \leq \bigg(\sum_{l=0}^{n}\norm{R_l}^2 + n+1\bigg)\bigg(\sum_{l=0}^{n}\frac{\norm{R_l}^2}{1+\norm{R_l}^2}\bigg) \\ & \leq 2(n+1) \sum_{l=0}^{n} \frac{\norm{R_l}^2}{1+\norm{R_l}^2} \\ & = (n+1)^2 \sum_{l=0}^{n}\frac{2\norm{R_l}^2}{(n+1)(1+\norm{R_l}^2)} \\ & = (n+1)^2 \biggl\{1- \frac{1}{n+1} \sum_{l=0}^{n}\frac{1-\norm{R_l}^2}{1+ \norm{R_l}^2}\biggr\}
\end{align*}
for any $n \in \mathbb{N}$. Further we observe that 
\begin{equation*}
	0 \leq \sum_{l=0}^{n}\frac{2\norm{R_l}^2}{(n+1)(1+\norm{R_l}^2)}  \leq \frac{1}{n+1}\sum_{l=0}^{n}\norm{R_l} \leq 1.
\end{equation*}
Therefore, the series $$q(n) = \frac{1}{n+1} \sum_{l=0}^{n}\frac{1-\norm{R_l}^2}{1+ \norm{R_l}^2}$$
must be convergent. If we proceed with exactly similar computations as in part (a) of the proof, then the inequality in part (b) can be easily obtained.
\end{pf}

\noindent\textbf{Acknowledgment:} 
The first author is supported by SERB-CRG, the second author is supported by CSIR (File No: 09/1059(0020)/2018-EMR-I), New Delhi, India, and the third author is supported by DST-INSPIRE Fellowship (IF 190721),  New Delhi, India.


\begin{thebibliography}{99}
	
	
	
	
	
	
	
	
	
	\bibitem{Ahamed-Allu-Halder-P3-2020} {\sc M. B. Ahamed, V. Allu} and {\sc H. Halder}, The Bohr Phenomenon for analytic functions on shifted disks, {\it  Ann. Fenn. Math.} {\bf 47} (2022), 103--120. 
	
	
	
	\bibitem{aizn-2000} {\sc L. Aizenberg}, Multidimensional analogues of Bohr's theorem on power series, \textit{Proc. Amer. Math. Soc.} {\bf 128} (2000), 1147--1155.
	
	
	\bibitem{aizenberg-2001} {\sc L. Aizenberg, A. Aytuna}  and {\sc P. Djakov}, Generalization of theorem on Bohr for bases in spaces of holomorphic functions of several complex variables, 
	{\it J. Math. Anal.Appl.} {\bf  258} (2001), 429--447.
	
	
	\bibitem{aizn-2007} {\sc L. Aizenberg}, Generalization of results about the Bohr radius for power series, {\it Stud. Math.}  {\bf 180}  (2007), 161--168.  
	
	
	
	
	
	
	
	
	
	
	\bibitem{Himadri-Vasu-P1} {\sc V. Allu} and {\sc H. Halder},  Bohr phenomenon for certain subclasses of Harmonic Mappings, {\it Bull. Sci. Math.} {\bf 173} (2021), 103053.
	
	\bibitem{Himadri-Vasu-P2} {\sc V. Allu} and {\sc H. Halder}, Bohr radius for certain classes of starlike and convex univalent functions, {\it J. Math. Anal. Appl.} {\bf 493} (2021), 124519.
	
	\bibitem{Himadri-Vasu-P9} {\sc V. Allu} and {\sc H. Halder}, Operator valued analogues of multidimensional Bohr's inequality, {\it Canadian Math. Bull.}  (2021), DOI: https://doi.org/10.4153/S0008439521001077.
	
	
	
	
	
	
	
	
	
	
	\bibitem{Ayt & Dja & BLMS & 2013} {\sc A. Aytuna} and {\sc P. Djakov}, Bohr property of bases in the space of entire functions and its generalizations, {\it Bull. London Math. Soc.} \textbf{45}(2)(2013), 411--420.
	
	\bibitem{bayart-advance-2014} {\sc F. Bayart, D. Pellegrino}, and {\sc J. B. Seoane-Sep\'{u}lveda}, The Bohr radius of the $n$-dimensional polydisk is equivalent to $\sqrt{(log \, n)/n}$, {\it Adv. Math.} {\bf 264} (2014), 726--746.
	
	
	\bibitem{bene-2004} {\sc C. B{\'e}n{\'e}teau}, {\sc A. Dahlner} and {\sc D. Khavinson}, Remarks on the Bohr phenomenon, {\it Comput. Methods Funct. Theory} \textbf{4}(1) (2004), 1-19.
	
	
	
	
		\bibitem{bhowmik-2018} {\sc B. Bhowmik} and {\sc N. Das},  Bohr phenomenon for subordinating families of certain univalent functions,  {\it J. Math. Anal. Appl.}  {\bf 462} (2018), 1087--1098.
	
	\bibitem{bhowmik-Edinburgh-2021} {\sc B. Bhowmik} and {\sc N. Das}, Bohr phenomenon for operator-valued functions, {\it Proc. Edinburgh Math. Soc.} {\bf 64} (2021), 72--86.
	
	\bibitem{bhowmik-Arch-Math-2021} {\sc B. Bhowmik} and {\sc N. Das}, A characterization of Banach spaces with nonzero Bohr radius, {\it 
		Arch. Math.} {\bf 116} (2021), 551--558.
	
		\bibitem{Blasco-Pavlovic-BLMS-2003} {\sc O. Blasco} and {\sc M. Pavlovic}, Complex convexity and vector-valued Littlewood-Paley inequalities, {\it Bull. London Math. Soc.} {\bf 35} (2003), 749--758.
		
	\bibitem{Blasco-OTAA-2010} {\sc O. Blasco}, The Bohr radius of a Banach space, {\it In Vector measures, integration and related topics, 5964, Oper. Theory Adv. Appl., 201, Birkh{\"a}user Verlag, Basel, 2010.}
	
		\bibitem{Blasco-Collect-2017} {\sc O. Blasco}, The $p$-Bohr radius of a Banach space, {\it Collect. Math.} {\bf 68} (2017), 87--100.
	
	
	\bibitem{boas-1997} {\sc H.P. Boas} and {\sc D. Khavinson}, Bohr's power series theorem in several variables, {\it Proc. Amer. Math. Soc.}  {\bf 125} (1997), 2975--2979.
	
	
	\bibitem{Bohr-1914} {\sc H. Bohr}, A theorem concerning power series,  {\it Proc. Lond. Math. Soc}. s2-13 (1914), 1--5.
	
	
	
	
	
	
	\bibitem{N-Das-JMAA-2022} {\sc N. Das}, Refinements of the Bohr and Rogosinski phenomena, {\it 
	J. Math. Anal. Appl.} {\bf 508} (2022), 125847, 10 pp.
	
	\bibitem{defant-2003} {\sc A. Defant, D. Garc\'{i}a}, and {\sc  M. Maestre}, Bohr power series theorem and local Banach space theory, {\it J. Reine Angew. Math.} {\bf 557} (2003), 173--197.
	
	\bibitem{defant-2006} {\sc A. Defant} and {\sc L. Frerick}, A logarithmic lower bound for multi-dimensional bohr radii, {\it Israel J. Math.} {\bf 152} (2006), 17--28.
	
	\bibitem{defant-2008} {\sc A. Defant, D. Garc\'{i}a, M. Maestre}, and {\sc D. P\'{e}rez-Garc\'{i}a} , Bohr's strip for vector-valued Dirichlet series, {\it Math. Ann.} {\bf 342} (2008), 533--555.
	
	
	\bibitem{defant-2012} {\sc A. Defant, M. Maestre}, and {\sc  U. Schwarting}, Bohr radii of vector valued holomorphic functions, {\it Adv. Math.} {\bf 231} (2012), 2837--2857.
	
	\bibitem{defant-2018} {\sc A. Defant, Mieczysław Mastyło}, and {\sc Antonio Pérez}, Bohr’s phenomenon for functions on the Boolean
	cube, {\it J. Funct. Anal.} {\bf 275} (2018), 3115--3147.
	
	\bibitem{dilworth-1986} {\sc S. Dilworth}, Complex convexity and geometry of Banach spaces, {\it Math. Proc. Camb. Phil. Soc.} {\bf 99} (1986), 495--506.

	
	\bibitem{Dixon & BLMS & 1995} {\sc P. G. Dixon}, Banach algebras satisfying the non-unital von Neumann inequality, {\it Bull. London Math. Soc.} \textbf{27} (1995), 359--362.
	
	\bibitem{Djakov & Ramanujan & J. Anal & 2000} {\sc P. B. Djakov} and {\sc M. S. Ramanujan}, A remark on Bohr's theorem and its generalizations, J. Anal. \textbf{8} (2000), 65--77.
	
	\bibitem{globevnik-1975} {\sc J. Globevnik}, On complex strict and uniform convexity, {\it Proc. Amer. Math. Soc.} {\bf 47} (1975), 175--178.
	
	
	
	
	
	
	
	
	
	
	
	
	
	
	
	
	
	
	
	
		\bibitem{Huang-Liu-Ponnu} {\sc Y. Huang, M-S. Liu}, and {\sc S. Ponnusamy}, Refined Bohr-type inequalities with area measure for bounded analytic functions, {\it Anal. Math. Phys.} {\bf 10} (2020), 50, 21 pp. 
	
	
	
	
	
	
	
	
	
		\bibitem{kayumov-2018-c} {\sc I. R. Kayumov} and {\sc S. Ponnusamy}, Improved version of Bohr's inequalities, 
		{\it C. R. Math. Acad. Sci. Paris} {\bf 358} (5) (2020), 615--620.
	
	
	
	
	
	
	
	
	
	
	
	
	
	
	
		\bibitem{Liu-2021} {\sc G. Liu}, {\sc Z. Liu} and {\sc S. Ponnusamy}, Refined Bohr inequality for bounded analytic functions, {\it Bull. Sci. Math.} {\bf 173} (2021), 103054, 20 pp.
	
	
	
	\bibitem{Liu-Pon-PAMS-2020} {\sc M. S. Liu} and {\sc S. Ponnusamy}, Multidimensional analogues of refined Bohr's inequality, {\it Proc. Amer. Math. Soc.} {\bf 149} (2021), 2133--2146.
	
	
	
	
	
	
	
	\bibitem{paulsen-2002} {\sc Vern I. Paulsen, Gelu Popescu} and {\sc Dinesh Singh}, On Bohr's inequality, {\it Proc. Lond. Math. Soc.} s3-85 (2002), 493--512.
	
	
	
	
	\bibitem{popescu-2019} {\sc G. Popescu}, Bohr inequalities for free holomorphic functions on polyballs, {\it Adv. Math.} {\bf 347} (2019), 1002-1053.
	
	
	
	
	
	
	\bibitem{rogosinski-1930} {\sc W. Rogosinski}, Uber Bildschranken bei Potenzreihen und ihren Abschnitten, {\it Math. Z.} {\bf 17} (1923), 260-276.
	
	
	
	
	
	
	
	
	
	
	
	\bibitem{thorp-PAMS-1967} {\sc E. Thorp} and {\sc R. Whitley}, The strong maximum modulus theorem for analytic functions into a Banach space, {\it Proc. Amer. Math. Soc.} {\bf 18} (1967), 640--646.
	
	
	
	
	
\end{thebibliography}
\end{document}